\documentclass[12pt,a4paper,final]{article}

\usepackage[latin1]{inputenc}

\usepackage[T1]{fontenc}
\usepackage[english,frenchb]{babel}
\usepackage{amsmath,amssymb,amscd,latexsym}

\usepackage[draft]{showkeys}

\usepackage{theorem}
\theorembodyfont{\slshape}
\newtheorem{Thm}{Th\'eor\`eme}[section]

\theoremstyle{plain}

\newtheorem{Def}[Thm]{D\'efinition}
\newtheorem{Prop}[Thm]{Proposition}
\newtheorem{Lem}[Thm]{Lemme}

\newtheorem{Ex}[Thm]{Exemple}
\newtheorem{Rem}[Thm]{Remarque}

\newenvironment{proof}{\noindent{\sc Preuve.}}{\hfill\qed}
\newenvironment{sproof}[1]{\noindent {\sc Preuve #1.}}{\hfill \qed}

\newcommand{\N}{\mathbb{N}}
\newcommand{\Z}{\mathbb{Z}}

\newcommand{\F}{\mathbb{F}}

\newcommand{\qed}{\hfill $\square$}
\newcommand{\Pres}[2]{\left\langle{#1}\ \big\vert\ {#2}\right\rangle}
\newcommand{\Engend}[1]{\left\langle {#1} \right\rangle}

\newcommand{\dom}{\operatorname{Dom}}

\title{MoyennabilitÈ intÈrieure et extensions HNN}
\author{Yves Stalder\footnote{FinancÈ par le Fonds national suisse de la recherche scientifique, subside n∞ 20-101469}}
\date{21 septembre 2005}

\begin{document}

\baselineskip=16pt

\maketitle

\selectlanguage{english}
\begin{abstract}
We present sufficient conditions for HNN extensions to be inner amenable, respectively ICC, which
give necessary and sufficient criteria among Baumslag-Solitar groups. We deduce that such a group, 
viewed as acting on its Bass-Serre tree, contains non trivial elements which fix unbounded subtrees.
\end{abstract}

\selectlanguage{french}
\begin{abstract}
On prÈsente des conditions suffisantes pour qu'une extension HNN soit intÈrieurement moyennable,
respectivement CCI, qui donnent des critËres nÈcessaires et suffisants parmi les groupes de Baumslag-Solitar. 
On en dÈduit qu'un tel groupe, vu comme groupe d'automorphismes de son arbre de Bass-Serre, possËde des 
ÈlÈments non triviaux qui fixent des sous-arbres non bornÈs.
\end{abstract}



\section*{Introduction}

…tant donnÈ un groupe $\Gamma$, l'Ètude de son algËbre de von Neumann $W^*(\Gamma)$ suggËre souvent des
propriÈtÈs et des questions intÈressantes concernant $\Gamma$. Par exemple, Murray et von Neumann avaient dÈj‡ remarquÈ
que $W^*(\Gamma)$ est un facteur de type $II_1$ si et seulement si le groupe $\Gamma$ est CCI, c'est-‡-dire \emph{‡ classes de
conjugaison infinies} ; ils avaient aussi dÈfini la <<~propriÈtÈ Gamma~>> pour un tel facteur. (Voir le lemme 5.3.4 et la
dÈfinition 6.1.1 de \cite{MvN}.) Plus tard, Effros \cite{Eff} a montrÈ que cette propriÈtÈ pour un facteur de type
$II_1$ de la forme $W^*(\Gamma)$ implique que le groupe $\Gamma$ est intÈrieurement moyennable (voir la section \ref{Prelim}
pour la dÈfinition ou [BH86] pour une introduction plus complËte ‡ cette notion).
Il rÈsulte immÈdiatement des dÈfinitions que les groupes non CCI et les groupes moyennables
sont intÈrieurement moyennables. Parmi les groupes CCI, il en existe qui sont
intÈrieurement moyennables, par exemple le cÈlËbre groupe $F$ de Thompson \cite{Jol}, et
d'autres qui ne le sont pas, par exemple les groupes libres non abÈliens.
Dans le cas des produits libres, l'Ètude de ces propriÈtÈs est particuliËrement simple.
En effet, si $\Gamma = H * K$ est produit libre de deux groupes $H$ et $K$ non rÈduits ‡ un
ÈlÈment, les trois propriÈtÈs suivantes sont Èquivalentes : (i) $\Gamma$ n'est pas CCI, (ii) $\Gamma$ est
intÈrieurement moyennable et (iii) $H$ et $K$ sont chacun d'ordre deux, de sorte que $\Gamma$ est
un groupe diÈdral infini. Pour des constructions plus gÈnÈrales, c'est un problËme ouvert
que de formuler (mÍme conjecturalement) des conditions nÈcessaires et suffisantes pour les
propriÈtÈs d'Ítre CCI et d'Ítre intÈrieurement moyennable. Notre objectif est d'Ètablir des
critËres s'appliquant au moins ‡ tous les groupes de Baumslag-Solitar (voir les exemples
2.4 et 3.2).

Pour un groupe $G$, nous posons $G^* = G\setminus \{1\}$ et nous notons $Z(G)$ le centre de $G$. Tous
les groupes considÈrÈs dans cet article sont supposÈs Ítre dÈnombrables.

\begin{Thm}\label{HNNCCI}
    Soit $\Gamma = HNN(\Lambda,H,K,\phi)$, avec $H \neq \Lambda$ ou $K \neq \Lambda$. Pour que
    $\Gamma$ soit ‡ classes de conjugaison infinies, il suffit que toute classe de conjugaison finie de
    $\Lambda$ incluse dans $(H\cap K)^*$ contienne un ÈlÈment qui n'est point fixe d'aucun homomorphisme $\phi^j$ 
    (avec $j\in\N^*$).

    En particulier, si l'un au moins des groupes $\Lambda$, $H$, $K$, $H\cap K$ est CCI ou si les homomorphismes $\phi^j$
    (avec $j\in\N^*$) sont sans point fixe non trivial, le groupe $\Gamma$ est ‡ classes de conjugaison infinies.
\end{Thm}
\begin{Prop}\label{HNNintMoy}
    Soit $\Gamma = HNN(\Lambda,H,K,\phi)$. Si pour tout $n\geqslant 1$ il existe $h_0^{(n)}, \ h_1^{(n)}, \ldots, h_n^{(n)}
    \in (Z(\Lambda) \cap H \cap K)^*$ tels que $h_i^{(n)} = \phi(h_{i-1}^{(n)})$
    pour $i = 1, \ldots, n$, alors $\Gamma$ est intÈrieurement moyennable.
\end{Prop}

On donne Ègalement un ÈnoncÈ que nous devons ‡ Guyan Robertson et qui assure la propriÈtÈ CCI pour
certains groupes agissant sur des arbres (voir la proposition \ref{CCIarbres}).

\medskip

De la moyennabilitÈ
intÈrieure, on peut encore dÈduire des propriÈtÈs de l'action de certaines extensions HNN sur leur arbre de
Bass-Serre:
\begin{Prop}\label{CorActionArbre}
    Soit $T$ un arbre et $\Gamma$ un sous-groupe de $\text{Aut}(T)$ contenant des automorphismes hyperboliques
    transverses. Si $\Gamma$ est intÈrieurement moyennable, il existe un ÈlÈment elliptique non trivial de $\Gamma$ qui
    fixe un sous-arbre non bornÈ de $T$.
\end{Prop}
Ce dernier ÈnoncÈ s'applique Ègalement ‡ la plupart des groupes de Baumslag-Solitar. Pour ces derniers, on
peut d'ailleurs exhiber explicitement des ÈlÈments qui fixent des sous-arbres non bornÈs.

\paragraph{Structure de l'article} La section \ref{Prelim} est consacrÈe aux rappels et dÈfinitions
nÈcessaires, la section \ref{CCI} ‡ la propriÈtÈ CCI et la section \ref{MoyInt} ‡ la moyennabilitÈ intÈrieure.
Ces deux derniËres sections sont agrÈmentÈes d'exemples. Enfin, la section \ref{MoyIntArbres} fait le lien
entre moyennabilitÈ intÈrieure et actions sur les arbres.

\paragraph{Remerciements} Je voudrais remercier Paul Jolissaint, qui m'a incitÈ ‡ Ètudier la moyennabilitÈ
intÈrieure des groupes de Baumslag-Solitar, ainsi que Pierre de la Harpe, dont les prÈcieuses
suggestions m'ont permis d'amÈliorer sensiblement la prÈsentation de cet article. Je remercie Ègalement Alain Valette, 
qui m'a montrÈ le lemme \ref{CCIlemArbres} et la proposition \ref{CCIarbres}, dus ‡ Guyan Robertson. Que les trois premiËres 
personnes soient en outre remerciÈes pour leurs commentaires sur les prÈcÈdentes versions de ce texte.

\section{PrÈliminaires}\label{Prelim}

\paragraph{MoyennabilitÈ intÈrieure et classes de conjugaison infinies} Soit $\Gamma$ un groupe (discret,
dÈnombrable). Le groupe $\Gamma$ agit sur $\Gamma^*$ par conjugaison, ce qui s'exprime par la formule
$\gamma\cdot x = \gamma x\gamma^{-1}$. Par abus de langage, on note encore $\gamma$ l'application $x\mapsto
\gamma \cdot x$.
\begin{Def}\label{CCIdef}
Un groupe $\Gamma$ est dit \emph{‡ classes de conjugaison infinies (CCI)} s'il est infini et si toutes les
orbites de l'action prÈcitÈe sont infinies, c'est-‡-dire si pour tout $x\in\Gamma^*$, l'ensemble
$\{ \gamma x\gamma^{-1} : \gamma\in\Gamma \}$ est infini.
\end{Def}
Une \emph{moyenne} sur $\Gamma$ est une forme linÈaire $m: \ell^\infty (\Gamma^*) \to \mathbb{C}$ positive
(c'est-‡-dire vÈrifiant $m(f) \geqslant 0$ si $f \geqslant 0$) et normalisÈe par $m(1)=1$.
\begin{Def}\label{MoyIntDef}
Un groupe $\Gamma$ est  dit \emph{intÈrieurement moyennable} si $\Gamma = \{1\}$ ou s'il existe une moyenne
sur $\Gamma$ qui soit invariante par automorphismes intÈrieurs, c'est-‡-dire telle que
$m(f \circ \gamma^{-1}) = m(f)$ pour tous $\gamma \in \Gamma$ et $f \in \ell^\infty (\Gamma^*)$.
\end{Def}

Un groupe non CCI est nÈcessairement intÈrieurement moyennable. En effet, si $X$ est une partie finie non vide de
$\Gamma^*$ stable par conjugaison, il suffit de poser
$$
m(f) = \frac{1}{|X|} \sum_{x\in X}f(x)
$$
pour obtenir une moyenne invariante par automorphismes intÈrieurs.

\paragraph{Extensions HNN} Soit $\Lambda$ un groupe, $H,K$ des sous-groupes et $\phi: H \to K$ un
isomorphisme. On dÈfinit l'\emph{extension HNN} de \emph{base} $\Lambda$ relativement
‡ $H$, $K$ et $\phi$ comme dans \cite[chapitre IV.2]{LynSch}:
$$
HNN(\Lambda,H,K,\phi) = \Pres{\Lambda, t}{t^{-1}ht = \phi(h) \ \forall h \in H} \ .
$$
Pour tout $j \geqslant 1$, dÈfinissons comme suit des homomorphismes $\phi^j$. On pose $\phi^1=\phi$
(en particulier $\dom (\phi^1) = \dom(\phi) = H$). Pour $j\geqslant 2$, on dÈfinit
rÈcursivement
$$
\dom (\phi^j) = \phi^{-1} \big( \dom(\phi^{j-1}) \cap K \big) \subseteq H  \ ,
$$
et naturellement $\phi^j(h) = \phi(\phi^{j-1}(h))$ pour tout $h\in \dom(\phi^j)$.

Soit $\gamma = \lambda_0 t^{\varepsilon_1} \lambda_1 \cdots t^{\varepsilon_n} \lambda_n$ un ÈlÈment de
$HNN(\Lambda,H,K,\phi)$, avec $n \geqslant 0$, $\lambda_i \in \Lambda$ et $\varepsilon_i = \pm 1$.
L'Ècriture est dite \emph{rÈduite} si elle ne contient aucun sous-mot de type $t^{-1} h t$ avec $h\in H$
ou $tkt^{-1}$ avec $k \in K$. L'ÈnoncÈ suivant est une variante du thÈorËme de la forme normale dans les extensions HNN
\cite[thÈorËme IV.2.1]{LynSch}. Nous nous y rÈfÈrerons (un peu abusivement) sous le nom de \emph{lemme de Britton}.
\begin{Lem}
Soit $\gamma = \lambda_0 t^{\varepsilon_1} \lambda_1 \cdots t^{\varepsilon_n} \lambda_n$ comme avant.
Si l'Ècriture est rÈduite et si $\gamma = 1$, alors $n=0$ et $\lambda_0$ est l'ÈlÈment neutre de $\Lambda$.
\end{Lem}
Gr‚ce au lemme IV.2.3 de \cite{LynSch}, on peut en outre dÈfinir une \emph{fonction longueur} $\ell :
HNN(\Lambda,H,K,\phi) \longrightarrow \N$ en posant $\ell(\gamma) = n$ si l'Ècriture $\gamma = \lambda_0 t^{\varepsilon_1}
\lambda_1 \cdots t^{\varepsilon_n} \lambda_n$ est rÈduite.

\paragraph{Groupes de Baumslag Solitar} Ils sont dÈfinis par les prÈsentations
$$
BS(m,n) = \Pres{a,b}{ab^ma^{-1} = b^n}
$$
o˘ $m$ et $n$ sont des entiers non nuls. Ce sont des extensions HNN du groupe $\Z$. Plus prÈcisÈment, on a
$BS(m,n) = HNN(\Z,n\Z,m\Z,\phi)$ o˘ $\phi: n\Z \to m\Z$ est dÈfini par $\phi(nk) = mk$. On vÈrifie facilement
que dans ce cas
$$
\dom(\phi^j) = n_1^j d \Z \ \text{ pour tout } j\geqslant 1 \ ,
$$
o˘ $d=\text{pgcd}(m,n)$ et $n_1 = \frac{n}{d}$.

\paragraph{Arbres et thÈorie de Bass-Serre} Pour les dÈfinitions des notions de graphe et d'arbre, on renvoie le lecteur
au chapitre I.2 de \cite{Ser}. Un automorphisme d'arbre $\alpha$ est dit \emph{sans inversion} si $\phi(e) \neq \bar e$
pour toute arÍte $e$. …tant donnÈ un tel automorphisme, la proposition 25 de \cite{Ser} montre que:
\begin{itemize}
  \item soit $\alpha$ est \emph{elliptique}, c'est-‡-dire fixe au moins un sommet de l'arbre;
  \item soit $\alpha$ est \emph{hyperbolique}, c'est-‡-dire tel qu'il existe une gÈodÈsique doublement infinie, appelÈe
  \emph{axe}, sur laquelle $\alpha$ induit une translation d'amplitude non nulle.
\end{itemize}
Deux automorphismes hyperboliques sont dits \emph{transverses} si leurs axes n'ont qu'un nombre fini (Èventuellement
nul) de sommets en commun.

\medskip

…tant donnÈ une extension $\Gamma = HNN(\Lambda,H,K,\phi)$, l'\emph{arbre de Bass-Serre} associÈ est le graphe donnÈ par
$$
Som(T) = \Gamma/\Lambda \ ; \ Ar(T) = \Gamma/H \sqcup \Gamma/K \ ; \ \overline{\gamma H} = \gamma t K \
; \ \overline{\gamma K} = \gamma t^{-1} H \ ;
$$
$$
(\gamma H)^- = \gamma \Lambda \ ; \ (\gamma H)^+ = \gamma t \Lambda \ ; \
(\gamma K)^- = \gamma \Lambda \ ; \ (\gamma K)^+ = \gamma t^{-1} \Lambda
$$
o˘, Ètant donnÈ une arÍte $e$, on note $e^-$ son origine et $e^+$ son sommet terminal.
Le chapitre I.5 de \cite{Ser} et le thÈorËme 12 en particulier assurent qu'il s'agit
bien d'un arbre.\footnote{Il convient de remarquer que la dÈfinition d'extension HNN dans \cite{Ser} ne
correspond pas ‡ la nÙtre qui est celle de \cite{LynSch} (on passe de l'une ‡ l'autre en remplaÁant $t$ par
$t^{-1}$).} On peut orienter $T$ de maniËre ‡ ce que l'action Èvidente de $\Gamma$ sur $T$ prÈserve l'orientation en posant
$$
Ar_+(T) = \Gamma/H \ ; \ (\gamma H)^- = \gamma \Lambda \ ; \ (\gamma H)^+ = \gamma t \Lambda \ .
$$
Relevons que cet arbre orientÈ est doublement rÈgulier: en chaque sommet  de $T$ les arÍtes sortantes sont
en bijection avec $\Lambda / H$  tandis que les arÍtes entrantes sont en bijection avec  $\Lambda / K$.

\section{Classes de conjugaison infinies}\label{CCI}

Le premier objectif de cette section est de prouver le thÈorËme \ref{HNNCCI}. On traite sÈparÈment
les ÈlÈments du sous-groupe $H\cap K$ et les autres.
\begin{Lem}\label{LemTypes}
    Soit $\Gamma = HNN(\Lambda,H,K,\phi)$, avec $H \cap K \neq \Lambda$. Pour tout \linebreak $\gamma \in \Gamma
    \setminus (H \cap K)$, la classe de conjugaison de $\gamma$ dans $\Gamma$ est infinie.
\end{Lem}
\begin{proof}
Si $\gamma\in \Lambda \setminus H$, les ÈlÈments $t^{-n}\gamma t^n$ (avec $n\in \N^*$) sont distincts deux ‡ deux par
le lemme de Britton. Si $\gamma\in \Lambda \setminus K$, il en est de mÍme des ÈlÈments $t^{n}\gamma t^{-n}$. Il ne reste
donc que le cas $\gamma \in \Gamma \setminus \Lambda$ ‡ traiter, pour lequel on supposera $K \neq \Lambda$, laissant
la vÈrification  du cas analogue $H \neq \Lambda$ au lecteur.

Soit $\gamma = \lambda_0 t^{\varepsilon_1} \lambda_1 \cdots t^{\varepsilon_{n}}\lambda_{n}$ une Ècriture
rÈduite (avec $n \geqslant 1$) et soit $C$ la classe de conjugaison de $\gamma$. Il suffit d'exhiber une
suite $(\gamma_j)$ d'ÈlÈments de $C$ telle que $\lim_{j\to\infty}\ell(\gamma_j) = +\infty$.

Si $\varepsilon_1 = -1$, on trouve $\sigma \in \Lambda$ tel que $\sigma \lambda_0 \not\in K$ gr‚ce ‡ l'hypothËse
$\Lambda \neq K$. On pose ensuite
$$
\gamma_1 = t^{-\varepsilon_n} \sigma \gamma \sigma^{-1} t^{\varepsilon_n} = t^{-\varepsilon_n} (\sigma\lambda_0)
t^{-1} \lambda_1 t^{\varepsilon_2} \lambda_2 \cdots t^{\varepsilon_{n}} (\lambda_{n} \sigma^{-1}) t^{\varepsilon_n}
$$
et la derniËre Ècriture est rÈduite par construction. Il suffit alors de poser $\gamma_j = t^{-\varepsilon_n}
\gamma_{j-1}t^{\varepsilon_n}$ pour $j\geqslant 2$ car il vient $\ell(\gamma_j)= n+2j \to \infty$ pour $j \to
\infty$.

Si $\varepsilon_n = 1$, on trouve $\tau \in \Lambda$ tel que $\lambda_n \tau \not\in K$ gr‚ce ‡ l'hypothËse
$\Lambda \neq K$. On pose ensuite
$$
\gamma_1 = t^{\varepsilon_1} \tau^{-1} \gamma \tau t^{-\varepsilon_1} = t^{\varepsilon_1} (\tau^{-1}\lambda_0)
t^{\varepsilon_1} \lambda_1 \cdots t^{\varepsilon_{n-1}} \lambda_{n-1} t (\lambda_{n} \tau) t^{-\varepsilon_1}
$$
et la derniËre Ècriture est rÈduite par construction. Il suffit alors de poser $\gamma_j = t^{\varepsilon_1}
\gamma_{j-1}t^{-\varepsilon_1}$ pour $j\geqslant 2$ car il vient $\ell(\gamma_j)= n+2j \to \infty$ pour $j \to
\infty$.

Enfin, si $\varepsilon_1 = 1$ et $\varepsilon_n = -1$, il suffit de poser $\gamma_j = t^j\gamma t^{-j}$ pour
obtenir $\ell(\gamma_j)= n+2j \to \infty$ pour $j \to \infty$.
\end{proof}

\begin{Lem}\label{LemTypeI}
Soit $\Gamma = HNN(\Lambda,H,K,\phi)$ et soit $\gamma\in (H \cap K)^*$. On note $C_\Gamma$ (respectivement $C_\Lambda$)
la classe de conjugaison de $\gamma$ dans $\Gamma$ (respectivement $\Lambda$).  Pour que la classe $C_\Gamma$ soit infinie,
il suffit qu'une au moins des conditions suivantes soit satisfaite:
\begin{enumerate}
  \item[(a)] la classe $C_\Lambda$ est infinie;
  \item[(b)] la classe $C_\Lambda$ contient un ÈlÈment qui n'est point fixe d'aucun homomorphisme $\phi^j$
  (avec $j\in \N^*$).
\end{enumerate}
\end{Lem}

On convient qu'un ÈlÈment situÈ hors du domaine de $\phi^j$ n'est pas un point fixe de cet homomorphisme.

\medskip

\begin{proof}
Si $C_\Lambda$ est infinie, il en est de mÍme de $C_\Gamma$ car $C_\Lambda \subseteq C_\Gamma$.

S'il existe un ÈlÈment $c$ de $C_\Lambda$ qui n'est point fixe d'aucun $\phi^j$ (pour $j\in \N^*$), on a
$t^{-n} c t^n \neq c$ pour tout $n \in \N^*$. En effet, la relation $t^{-n} c t^n = c$ impliquerait
$\phi^{n}(c) = c$ gr‚ce au lemme de Britton. Par consÈquent, l'ensemble $\{t^{-n} c t^n : n\in \N\}$
est infini. Comme il est inclus dans la classe $C_\Gamma$, cette derniËre est Ègalement infinie.
\end{proof}

\medskip

\begin{sproof}{du thÈorËme \ref{HNNCCI}}
Soit $\gamma \in \Gamma^*$. Si la classe de conjugaison de $\gamma$ est contenue dans $(H\cap K)^*$ les hypothËses
impliquent qu'une des conditions (a),(b) du lemme \ref{LemTypeI} est satisfaite, ce qui entraÓne que la classe de conjugaison
de $\gamma$ est infinie. Dans le cas contraire, on obtient la mÍme conclusion gr‚ce au lemme \ref{LemTypes}.
\end{sproof}

\medskip

La proposition suivante s'applique ‡ un des cas d'extension HNN non couverts par le thÈorËme \ref{HNNCCI}, c'est-‡-dire au cas
d'un produit semi-direct $\Z\ltimes_{\phi} \Lambda$ associÈ ‡ un automorphisme $\phi$ de $\Lambda$ (cas o˘
$H=\Lambda = K$).

\begin{Prop}\label{HNNCCI2}
Soit $\Gamma = HNN(\Lambda,H,K,\phi)$. Si $\Lambda$ est non trivial et si les homomorphismes $\phi^j$
(avec $j\in\N^*$) sont sans point fixe non trivial, alors $\Gamma$ est CCI.
\end{Prop}
\begin{proof}
Par le thÈorËme \ref{HNNCCI}, il ne reste que le cas $H=\Lambda=K$ ‡ traiter. Tout ÈlÈment de $\Gamma$ peut
alors s'Ècrire sous la forme $\lambda t^n$ avec $\lambda \in \Lambda$ et $n\in\Z$. Supposons par
l'absurde qu'il existe un ÈlÈment $\gamma \in \Gamma^*$ dont la classe de conjugaison est finie.

On peut supposer que $\gamma = \lambda t^n$ avec $\lambda \in \Lambda^*$. En effet, si $\gamma = t^n$ avec
$n\in\Z^*$, on prend $\mu \in \Lambda^*$ (ce qui est licite car $\Lambda$ est non trivial) et on n'a plus
qu'‡ remplacer $\gamma$ par $\mu^{-1} \gamma \mu = \mu^{-1} \phi^{-n}(\mu) t^n$ puisque $\phi^{-n}$ est sans
point fixe non trivial. (Remarquons que la dÈfinition des itÈrÈs de $\phi$ ne pose pas de problËme ici puisqu'il
s'agit d'un automorphisme de $\Lambda$.)

Comme la classe de conjugaison de $\gamma$ est finie, il en est de mÍme de l'ensemble
$
\{t^{-j} \gamma t^j : j\in \N \} = \{\phi^j(\lambda) t^n : j \in \N \}
$,
ce qui montre qu'il existe $j \in \N^*$ tel que $\phi^j(\lambda) = \lambda$. On a obtenu la contradiction
cherchÈe (avec les hypothËses).
\end{proof}

\medskip

On passe maintenant ‡ des exemples d'extensions HNN ‡ classes de conjugaison infinies. Bien que nous n'en
ayons pas trouvÈ trace dans la littÈrature, ils sont peut-Ítre dÈj‡ connu des experts.

\begin{Ex}
…tant donnÈs deux entiers non nuls $m$ et $n$, le groupe de Baumslag-Solitar $BS(m,n)$ est CCI si et
seulement si $|m| \neq |n|$.
\end{Ex}
\begin{proof}
Si $m=n$, l'ÈlÈment $b^m$ est central dans $BS(m,m)$ et donc seul dans sa classe de conjugaison.
Si $m=-n$, l'ensemble $\{ b^m, b^{-m} \}$ est une classe de conjugaison finie.

Si $|m| \neq |n|$, le thÈorËme \ref{HNNCCI} implique que $BS(m,n)$ est CCI.
En effet, les homomorphismes $\phi^j$ associÈs ‡ la structure d'extension HNN de $BS(m,n)$ sont sans point fixe non trivial.
\end{proof}

\begin{Ex}
Si $\Gamma = HNN(\Lambda,H,K,\phi)$ avec $\Lambda = \Z^d = H$, les conditions suivantes sont Èquivalentes:
\begin{enumerate}
  \item[(i)] le groupe $\Gamma$ est CCI;
  \item[(ii)] aucune des valeurs propres complexes de l'endomorphisme $\phi$ n'est une racine de l'unitÈ.
\end{enumerate}
\end{Ex}
\begin{proof}
$(i) \Rightarrow (ii)$: Par contraposition, supposons que $\phi$ possËde une valeur propre $\mu$ telle que
$\mu^j = 1$ ($j\in\N^*$). Comme $1$ est valeur propre de $\phi^j$, celui-ci possËde un point fixe non nul dans
$\mathbb{Q}^d$. Il possËde donc aussi un point fixe non nul dans $\Z^d$, qu'on note $\lambda$. Par suite, l'ensemble $\{\lambda,
\phi(\lambda), \ldots, \phi^{j-1}(\lambda)\}$ est une classe de conjugaison finie de $\Gamma$.

$(ii) \Rightarrow (i)$: Pour tout $j\in \N^*$, l'endomorphisme $\phi^j$ est sans point fixe non trivial car $1$ n'en est
pas valeur propre. La proposition \ref{HNNCCI2} assure donc que $\Gamma$ est CCI.
\end{proof}

\medskip

Nous donnons maintenant un dernier exemple, qui montre que la condition du thÈorËme \ref{HNNCCI} n'est pas nÈcessaire.

\begin{Ex}
Posons $\Lambda = BS(m,m)$, avec $m\in \Z, |m| \geqslant 2$. Les ÈlÈments $a$ et $b^m$ engendrent un sous-groupe
isomorphe ‡ $\Z^2$, qu'on note $H$. On considËre $\Gamma = HNN(\Lambda,H,H,\phi)$, o˘ $\phi: H \to H$ est dÈfini
par $\phi(a) = b^m$ et $\phi(b^m) = a$.
\begin{enumerate}
  \item[(1)] Le groupe $\Gamma$ est CCI.
  \item[(2)] Il existe une classe de conjugaison finie de $\Lambda$ contenue dans $H^*$ et dont tous les ÈlÈments sont 
  des points fixes d'homomorphismes de la forme $\phi^j$, avec $j\in \N$.
\end{enumerate}
\end{Ex}
\begin{proof}
(1) Soit $\gamma$ un ÈlÈment de $\Gamma$ dont la classe de conjugaison est finie. On doit montrer que $\gamma
= 1$. Par le lemme \ref{LemTypes}, $\gamma$ est un ÈlÈment de $H$, et donc Ègalement de $\Lambda$. Mais
on peut Ècrire $\Lambda = HNN(\Z,m\Z,m\Z,\text{id})$, o˘ la base $\Z$ correspond au sous-groupe engendrÈ par $b$.
¿ nouveau gr‚ce au lemme \ref{LemTypes} on a $\gamma= b^{\ell m}$ avec $\ell \in \Z$. Par suite,
$t^{-1}\gamma t = a^\ell$, ce qui montre que la classe de conjugaison de $a^\ell$ est finie. Toujours gr‚ce au lemme
\ref{LemTypes}, on obtient alors
$\ell = 0$, ce qui donne finalement $\gamma = 1$ comme dÈsirÈ.

(2) L'ÈlÈment $b^m$ de $H^*$ est central dans $\Lambda$. Il est donc seul dans sa classe de conjugaison. De plus,
on a $\phi^2(b^m) = b^m$.
\end{proof}

\medskip

Les deux ÈnoncÈs suivants, dus ‡ Guyan Robertson, montrent plus gÈnÈralement que certains groupes
d'automorphismes d'arbres sont ‡ classes de conjugaison infinies.
\begin{Lem}\label{CCIlemArbres}
Soit $\Gamma$ un groupe agissant fidËlement et isomÈtriquement sur un espace mÈtrique $(X,d)$ de telle maniËre
que le quotient $\Gamma \backslash X$ soit fini. Si $\gamma$ est un ÈlÈment de $\Gamma$ dont la classe de
conjugaison est finie, alors il existe $K_\gamma >0$ tel que $d(x,\gamma x) \leqslant K_\gamma$ pour tout $x \in X$.
\end{Lem}
\begin{proof}
Soit $F\subseteq X$ un ensemble de reprÈsentants des orbites; cet ensemble est fini. La classe de conjugaison
de $\gamma$ Ètant finie, on peut poser
$$
K_\gamma = \max \{ d(y,g\gamma g^{-1}y) : g\in \Gamma, y\in F \} \ .
$$

Si $x$ est un ÈlÈment quelconque de $X$, on trouve $g$ dans $\Gamma$ et $y$ dans $F$ tels que $y=gx$. Alors
$d(x,\gamma x) = d(g^{-1}y, \gamma g^{-1}y) = d(y,g\gamma g^{-1}y) \leqslant K_\gamma$.
\end{proof}
\begin{Prop}\label{CCIarbres}
Soit $X$ un arbre dont tous les sommets ont au moins trois voisins. Si un groupe $\Gamma$ agit fidËlement
sur $X$ de telle maniËre que le quotient $\Gamma \backslash X$ soit fini, alors il est CCI.
\end{Prop}
\begin{proof}
Soit $\gamma$ un ÈlÈment de $\Gamma$ dont la classe de conjugaison est finie. On veut montrer que
$\gamma = 1$. Comme l'action est fidËle, il suffit de montrer que $\gamma$ agit trivialement sur $X$.

On commence par supposer (par l'absurde) que l'ÈlÈment $\gamma$ est hyperbolique. Comme tout sommet de $X$ a au moins
trois voisins, on peut trouver un sommet $x_n$ ‡ distance $n$ de l'axe pour tout $n \geqslant 1$. Mais alors,
$d(x_n, \gamma x_n) \geqslant 2n$, ce qui contredit le lemme \ref{CCIlemArbres}.

Par consÈquent, $\gamma$ fixe au moins un sommet de $X$. Notons $T$ le sous-arbre de $X$ formÈ des points
fixes de $\gamma$. On suppose par l'absurde que $T\neq X$. Dans ce cas, comme tout sommet est de degrÈ
au moins trois, on trouve $x_n$ ‡ distance $n$ de $T$ pour tout $n \geqslant 1$. Comme $\gamma$ ne fixe aucun
sommet de la gÈodÈsique joignant $x_n$ ‡ $T$, on a $d(x_n,\gamma x_n) = 2n$ pour tout $n$, ce qui contredit le lemme
\ref{CCIlemArbres}.

Ainsi, $T=X$, ce qui signifie que $\gamma$ agit trivialement sur $X$.
\end{proof}

\medskip

La proposition \ref{CCIarbres} s'applique en particulier au cas d'une extension $\Gamma = HNN(\Lambda,H,K,\phi)$ telle
que $[\Lambda:H] \geqslant 3$, $[\Lambda:K] \geqslant 3$ et $\bigcap_{\gamma\in\Gamma} \gamma^{-1}H\gamma =
\{1\}$, et de mÍme au cas d'un produit libre $H*_L K$ avec amalgamation au-dessus d'un sous-groupe commun $L$ de $H$ et
$K$,
tel que $[H:L] \geqslant 3$, $[K:L] \geqslant 3$ et $\bigcap_{\gamma\in\Gamma} \gamma^{-1}L\gamma = \{1\}$.
Terminons par un exemple d'extension HNN montrant que le thÈorËme \ref{HNNCCI} n'est pas consÈquence
de la proposition \ref{CCIarbres}.
\begin{Ex}
ConsidÈrons le groupe $\Lambda = \F_2 \times (\oplus_{\Z}\Z)$, de prÈsentation
$$
\Lambda = \Pres{a,b,e_i \ (i\in\Z)}{[a,e_i],[b,e_i],[e_i,e_j] \ (i,j\in\Z)} \ .
$$
Soit $H$ le sous-groupe engendrÈ par les $e_i$ et $\phi:H \to H$ l'homomorphisme dÈfini par la formule
$\phi(e_i) = e_{i+1}$. Soit encore $\Gamma = HNN(\Lambda,H,H,\phi)$.

Alors les homomorphismes $\phi^j$
(avec $j\in\N^*)$ sont sans point fixe non trivial (en particulier $\Gamma$ est CCI par le thÈorËme \ref{HNNCCI}),
mais l'action de $\Gamma$ sur son arbre de Bass-Serre n'est pas fidËle.
\end{Ex}
\begin{proof}
Il est Èvident que les homomorphismes $\phi^j$ (avec $j\in\N^*)$ sont sans point fixe non trivial. D'autre part,
le sous-groupe $H$ est normal dans $\Gamma$ par construction. Notons $T$ l'arbre de Bass-Serre de $\Gamma$.
On a $Ar_+(T) = \Gamma / H$, de sorte que  le sous-groupe $H$ fixe toutes les arÍtes de l'arbre de
Bass-Serre. Ainsi, l'action de $\Gamma$ sur $T$ n'est pas fidËle.
\end{proof}

\section{MoyennabilitÈ intÈrieure}\label{MoyInt}

Nous nous intÈressons maintenant ‡ montrer que certaines extensions HNN sont intÈrieurement moyennables.
La preuve qui suit s'inspire du corollaire 2.5. de \cite{Jol}.

\medskip

\begin{sproof}{de la proposition \ref{HNNintMoy}} Sans nuire ‡ la gÈnÈralitÈ, on peut supposer que le groupe
$\Gamma$ est CCI. Il s'ensuit que pour tout $n$, les ÈlÈments $h_0^{(n)}, \ldots, h_n^{(n)}$ sont deux ‡ deux
distincts. En effet, si $h_i^{(n)} = h_j^{(n)}$ (avec $i<j$), l'ensemble $\{ h_i^{(n)}, \ldots h_{j-1}^{(n)}\}$
est une classe de conjugaison finie.

Nous allons construire une moyenne sur $\Gamma$
qui soit invariante par automorphismes intÈrieurs. Soit $\omega$ un ultrafiltre sur $\N^*$ plus fin
que le filtre de FrÈchet\footnote{Il est Èquivalent de demander que $\omega$ soit un ultrafiltre libre.}.
On pose pour $f \in \ell^\infty (\Gamma^*)$ :
\begin{eqnarray*}
    \mu_n(f)      & = & \frac{1}{n} \sum\limits_{i=0}^{n-1} f(h_i^{(n)}) \quad (n \in \N^*)   \\
    \mu_\omega(f) & = & \underset{n\to\omega}{\lim} \ \mu_n(f) \ .
\end{eqnarray*}
La limite existe car la suite $(\mu_n(f))_n$ est bornÈe et $\omega$ est un ultrafiltre. La forme linÈaire
$\mu_\omega$ sur $\ell^\infty(\Gamma^*)$ est positive et normalisÈe. Il reste ‡ voir
qu'elle est invariante par automorphismes intÈrieurs. Comme $\Gamma$ est engendrÈ par $\Lambda$ et $t$, cela revient ‡
montrer:
\begin{enumerate}
  \item[(1)] $\mu_\omega(f \circ S_\lambda) = \mu_\omega(f)$ pour tout $\lambda \in \Lambda$, o˘ $S_\lambda(\gamma)
  = \lambda^{-1} \gamma \lambda$;
  \item[(2)] $\mu_\omega(f \circ T) = \mu_\omega(f)$, o˘ $T(\gamma) = t^{-1} \gamma t$.
\end{enumerate}
Comme les $h_i^{(n)}$ sont centraux dans $\Lambda$, on a $S_\lambda(h_i^{(n)}) = h_i^{(n)}$ pour tout $n \in \N^*$ et
pour tout $i= 1,\ldots,n$. Cela montre (1). Pour montrer (2), on constate d'abord que:
\begin{eqnarray*}
    \mu_n(f\circ T) & = & \frac{1}{n} \sum\limits_{i=0}^{n-1} f\circ T(h_i^{(n)})
                        = \frac{1}{n} \sum\limits_{i=0}^{n-1} f(t^{-1} h_i^{(n)} t)   \\
                    & = & \frac{1}{n} \sum\limits_{i=0}^{n-1} f\circ \phi(h_i^{(n)})
                        = \frac{1}{n} \sum\limits_{i=0}^{n-1} f(h_{i+1}^{(n)})          \\
                    & = & \frac{1}{n} \sum\limits_{i=1}^n f(h_i^{(n)})
                        = \mu_n(f) + \frac{1}{n} \big( f(h_n^{(n)}) - f(h_0^{(n)}) \big)
\end{eqnarray*}
Le dernier terme tend vers $0$ selon $\omega$. Il suit
que $\mu_\omega(f \circ T) = \mu_\omega(f)$.
\end{sproof}

\begin{Rem}
    Pour prouver la proposition \ref{HNNintMoy}, il est plus rapide de vÈrifier la condition (F) du thÈorËme
    1 de \cite{BedHar}. On suppose comme avant que $\Gamma$ est CCI, de sorte que les $h_i^{(n)}$ sont deux ‡
    deux distincts. On constate alors qu'on peut prendre $F_n = \{ h_1^{(n)}, \ldots, h_{n-1}^{(n)} \}$ pour
    $n \geqslant 2$. En effet $|\lambda F_n \lambda^{-1} \vartriangle F_n| = 0$ pour $\lambda \in \Lambda$ et
    $|t^{\pm 1} F_n t^{\mp 1} \vartriangle F_n| = 2$. Par suite, $|\gamma F_n \gamma^{-1} \vartriangle F_n|
    \leqslant 2\ell(\gamma)$ pour tous $\gamma \in \Gamma$ et $n\geqslant 2$. Pour tout $\gamma\in \Gamma$ il
    vient donc
    $$
    \underset{n \to \infty}{\text{lim}} \frac{|\gamma F_n \gamma^{-1} \vartriangle F_n|}{|F_n|} = 0 \ .
    $$
\end{Rem}

L'exemple suivant infirme la proposition 4.3 de \cite{BegCec}.\footnote{Pour dÈmontrer leur ÈnoncÈ, les
auteurs ont voulu utiliser la proposition 7 de \cite{BedHar}.
Dans la vÈrification des hypothËses, ils ont supposÈ ‡ tort que les ÈlÈments du groupe dÈfinissant des
automorphismes elliptiques de l'arbre de Bass-Serre fixaient un unique sommet. Voir aussi l'exemple \ref{BSfix}}
\begin{Ex}\label{BSintMoy}
Le groupe $BS(m,n)$ est intÈrieurement moyennable quels que soient les entiers non nuls $m$ et $n$.
\end{Ex}
\begin{proof}
Le groupe $BS(m,n)$ Ètant une extension HNN, il suffit de vÈrifier qu'on peut appliquer la proposition
\ref{HNNintMoy}. Pour tout $k \geqslant 1$, il suffit de poser $h_i^{(k)} =
b^{m^{i+1} n^{k-i+1}}$ pour obtenir $h_i^{(k)} = \phi(h_{i-1}^{(k)})$ pour tout $i = 1, \ldots, k$.
En outre, la base Ètant commutative, ces ÈlÈments sont bien centraux.
\end{proof}

\begin{Ex}\label{AscintMoy}
Si $\Lambda$ est abÈlien et si $H=\Lambda$, alors $\Gamma = HNN(\Lambda,H,K,\phi)$ est intÈrieurement
moyennable.
\end{Ex}
\begin{proof}
Si $\Lambda$ est trivial, on a $\Gamma = \Z$ et ce groupe est intÈrieurement moyennable. Sinon, il suffit de
prendre $\lambda \in \Lambda^*$ et de poser $h_i^{(n)} = \phi^i(\lambda)$ pour pouvoir appliquer la
proposition \ref{HNNintMoy}
\end{proof}

\medskip

Nous terminons cette section en donnant des exemples qui infirment les points (d) et (e) du thÈorËme 5 de
\cite{BedHar}\footnote{Dans cet article Ègalement, les auteurs ont supposÈ ‡ tort que les ÈlÈments elliptiques
(pour l'action sur l'arbre de Bass-Serre associÈ) fixaient un seul sommet. Dans le cas (c) du thÈorËme visÈ,
par contre, cette hypothËse est vÈrifiÈe car les stabilisateurs des arÍtes sont triviaux.}, qui donnaient des
obstructions ‡ la moyennabilitÈ intÈrieure pour des groupes du type $\Gamma = H *_A K$, respectivement $\Gamma =
HNN(H,A,B,\phi)$.

De nouvelles obstructions ‡ la moyennabilitÈ intÈrieure pour les extensions HNN et produits amalgamÈs
peuvent par contre Ítre obtenues en appliquant la proposition \ref{CorActionArbre} aux actions sur les
arbres de Bass-Serre correspondants.

\begin{Rem}
L'hypothËse principale des points (d) et (e) dans \cite{BedHar} Ètait:  <<~pour toute partie finie
$F\subseteq \Gamma^*$, il existe $\gamma \in \Gamma$ tel que $\gamma F \gamma^{-1} \cap A = \varnothing$~>>.
GÈomÈtriquement parlant, cette condition assure que les actions des groupes sur leur arbre de Bass-Serre
et sur le bord de ce dernier sont fortement fidËles. Voir par exemple le lemme 9, sa preuve et la remarque
prÈcÈdent la proposition 11 dans \cite{Har}.
\end{Rem}

Le groupe $BS(2,3)$ fournit un contre-exemple au point (e). En effet, on a vu dans l'exemple \ref{BSintMoy}
que ce groupe est intÈrieurement moyennable tandis que le lemme qui suit montre que les hypothËses de (e)
sont satisfaites.
\begin{Lem}\label{HNNFidFort}
    Si $\Gamma = HNN(\Lambda,H,K,\phi)$ est tel que pour tout $\lambda \in \Lambda$, il existe $k \in \N$
    avec $\phi^k(\lambda) \not \in H$, alors pour toute partie finie $F$
    de $\Gamma^*$ et pour $n$ suffisamment grand, on a $t^{-n} F t^n \cap \Lambda = \varnothing$.
\end{Lem}

\begin{proof} Il suffit de montrer que pour tout $\gamma \in \Gamma^*$ et pour $n$
suffisamment grand, on a $t^{-n} \gamma t^n \not \in \Lambda$. On prend donc $\gamma \in \Gamma^*$ et on suppose
pour Èviter la trivialitÈ qu'il existe $n_0 \in \N$ avec $t^{-n_0} \gamma t^{n_0} \in \Lambda$. Pour $n \in \N$,
on pose $\gamma_n = t^{-n_0 - n} \gamma t^{n_0 + n}$. Par hypothËse, il existe $k \in \N$ tel que $\gamma_{k}
\in \Lambda \setminus H$. Les Ècritures $\gamma_{k+n} = t^{-n} \gamma_k t^n$ sont alors rÈduites pour tout
$n \geqslant 1$. Donc pour $n > n_0 + k$, on a bien $t^{-n} \gamma t^n \not \in \Lambda$.
\end{proof}

\medskip

Le contre-exemple au point (d) est le suivant:
\begin{Ex}\label{contrex5d}
ConsidÈrons le groupe
$$
\Gamma = \Pres{a_1,a_2,b}{a_1 b^2 a_1^{-1} = b^3, a_2 b^2 a_2^{-1} = b^3} = BS(2,3)*_\Z BS(2,3) \ :
$$
\begin{enumerate}
 \item[(1)] il est intÈrieurement moyennable;
 \item[(2)] pour tout ensemble fini $F$ d'ÈlÈments de $\Gamma^*$, il existe $\gamma \in \Gamma$ tel que $\gamma F
 \gamma^{-1} \cap \langle b \rangle = \varnothing$.
\end{enumerate}
\end{Ex}

\begin{proof} (1) Il est facile de voir que la suite de parties
$$
F_n = \{ b^{2^n}, b^{2^{n-1}\cdot 3}, \ldots, b^{2\cdot 3^{n-1}}, b^{3^n} \}
$$
satisfait la condition (F) du thÈorËme 1 de \cite{BedHar}.

(2) On se donne $g \in \Gamma^*$. Il suffit de montrer que, pour $n$ suffisamment grand, l'ÈlÈment $(a_1a_2)^n
g (a_1a_2)^{-n}$ n'est pas une puissance de $b$. Pour Ècarter le cas trivial, on peut supposer qu'il existe
$n_0$ tel que $(a_1a_2)^{n_0} g (a_1a_2)^{-n_0}$ soit une puissance de $b$. Vu les relations qui dÈfinissent
$\Gamma$ il existe $n_1 \geqslant n_0$ tel que $(a_1a_2)^{n_1} g (a_1a_2)^{-n_1} = b^k$ avec $k$ non
divisible par $4$. On aura dËs lors:
\begin{eqnarray*}
(a_1a_2)^{n_1 + 1} g (a_1a_2)^{-n_1 -1} & = & a_1 b^{3\frac{k}{2}}a_1^{-1} \hspace{1cm} \text{ si } k \text{ est pair;} \\
(a_1a_2)^{n_1 + 1} g (a_1a_2)^{-n_1 -1} & = & a_1a_2 b^k a_2^{-1}a_1^{-1} \hspace{2mm} \text{ si } k \text{ est impair.}
\end{eqnarray*}
Les expressions $a_1 b^{3\frac{k}{2}}a_1^{-1}$ et $a_2 b^k a_2^{-1}$ ne sont pas rÈductibles dans $BS(2,3)$
(par le lemme de Britton). Pour $n > n_1 + 1$, l'ÈlÈment $(a_1a_2)^{n} g (a_1a_2)^{-n}$ n'est donc pas
une puissance de $b$ \cite[thÈorËme IV.2.6]{LynSch}.
\end{proof}

\section{MoyennabilitÈ intÈrieure et arbres}\label{MoyIntArbres}

…tant donnÈ un arbre $T$, on note $\partial T$ l'espace des bouts de $T$. Remarquons que les automorphismes
de $T$ dÈfinissent des homÈomorphismes de l'espace $T\cup \partial T$. De plus un automorphisme de $T$ est
hyperbolique si et seulement s'il dÈfinit un homÈomorphisme hyperbolique de $T\cup \partial T$ au sens de
\cite[pp.~151--152]{BedHar}.
Notons Ègalement que deux automorphismes hyperboliques de $T$ sont transverses au sens de la section \ref{Prelim} si et
seulement si les homÈomorphismes hyperboliques de $T\cup \partial T$ associÈs sont transverses au sens de \cite{BedHar}.
Le lien entre homÈomorphismes hyperboliques et moyennabilitÈ intÈrieure est le suivant:
\begin{Prop}\label{BH}\cite[proposition 7]{BedHar}
Soit $\Omega$ un espace topologique sÈparÈ infini et $\Gamma$ un groupe d'homÈomorphismes de $\Omega$.
Si $\Gamma$ contient des homÈomorphismes hyperboliques transverses de $\Omega$ et s'il existe une application
$\Gamma$-Èquivariante $\delta : \Gamma^* \to \Omega$ (o˘ $\Gamma$ opËre sur $\Gamma^*$ par conjugaison),
alors $\Gamma$ n'est pas intÈrieurement moyennable.
\end{Prop}

Rappelons en outre que les arbres sont des espaces mÈtriques complets satisfaisant l'inÈgalitÈ de la mÈdiane
\cite[chapitre 3.b]{HaVa}. Le lemme 3.8 de \cite{HaVa} assure que toute partie bornÈe et non
vide d'un arbre possËde un \emph{centre}, qui est le centre de la boule fermÈe de rayon minimal qui la
contienne.

\medskip

\begin{sproof}{de la proposition \ref{CorActionArbre}} Quitte ‡ remplacer $T$ par sa subdivision
barycentrique, on peut supposer que $\Gamma$ agit sans inversion.

Par contraposition, on suppose que, pour tout ÈlÈment
elliptique $\gamma$ de $\Gamma$, le sous-arbre $T_{\gamma}$ formÈ des points fixes de $\gamma$ est bornÈ. L'existence
d'ÈlÈments hyperboliques implique que $T$ n'est pas rÈduit ‡ un sommet. Ainsi $\Omega = T \cup \partial T$ est un espace
topologique sÈparÈ infini. ConsidÈrons l'application
$
\delta : \Gamma^* \longrightarrow \Omega
$,
dÈfinie comme suit:
\begin{itemize}
    \item Si $\gamma$ est hyperbolique, $\delta(\gamma)$ est le point attractif de $\gamma$;
    \item Si $\gamma$ est elliptique,  $\delta(\gamma)$ est le centre de $T_{\gamma}$.
\end{itemize}
Si on fait agir $\Gamma$ sur $\Gamma^*$ par automorphismes intÈrieurs, l'application $\delta$ est
$\Gamma$-Èquivariante. Gr‚ce ‡ la proposition \ref{BH}, on en dÈduit que $\Gamma$ n'est pas
intÈrieurement moyennable, ce qui contredit les hypothËses.
\end{sproof}

\begin{Ex}\label{BSfix}
Pour tous $m,n \in \Z^*$, il existe $\gamma \in BS(m,n)$ tel que $\gamma$ fixe un sous-arbre
non bornÈ de l'arbre de Bass-Serre associÈ.
\end{Ex}
\begin{proof}
Remarquons pour commencer que, dans le cas qui nous occupe, l'arbre de Bass-Serre est localement fini. Les sous-arbres
non bornÈs coÔncident donc avec les sous-arbres infinis. On distingue trois cas.

(1) Si $m = kn$ avec $k\in\Z^*$, on peut prendre
$$
\gamma = b^n = a b^m a^{-1} = a b^{kn} a^{-1} = a^2 b^{km} a^{-2} = a^2 b^{k^2 n} a^{-2} = \ldots
$$
On a $\gamma = a^\ell b^{k^\ell n} a^{-\ell} \ \forall \ell \in \N$. L'ÈlÈment $\gamma$ fixe (au moins) les sommets $a^\ell
\Engend{b}$ o˘ $\ell$ parcourt $\N$ et donc un sous-arbre non bornÈ de $T$.

(2) Si $n = km$ on constate de maniËre analogue que $\gamma = b^m$ fixe les sommets $a^{-\ell} \Engend{b}$
o˘ $\ell$ parcourt $\N$.

(3) Si les nombres $m$ et $n$ ne sont pas multiples l'un de l'autre (en particulier on a $|m|,|n| \geqslant 2$),
on prend $\gamma = b^n$. Remarquons que
$$
(aba^{-1}b^{-1})^\ell \gamma (aba^{-1}b^{-1})^{-\ell} = \gamma \text{ pour tout } \ell \in\N \ .
$$
Par suite, $\gamma$ fixe tous les sommets $(aba^{-1}b^{-1})^{k} \Engend{b}$.
Les expressions $(aba^{-1}b^{-1})^{k}$ Ètant rÈduites, le lemme de Britton implique que ces sommets dont distincts
deux ‡ deux. L'ÈlÈment $\gamma$ fixe donc un sous-arbre non bornÈ.
\end{proof}
\begin{Rem}
On laisse le soin au lecteur de vÈrifier par lui-mÍme que la proposition \ref{CorActionArbre} s'applique
si $|m|,|n| \geqslant 2$.
\end{Rem}

\bibliographystyle{alpha}

\medskip

Adresse de l'auteur:

\vspace{2mm}

Institut de Math\'ematiques - Universit\'e de Neuch\^atel

Rue Emile Argand 11

CH-2007 Neuch\^atel - SUISSE

\vspace{2mm}

yves.stalder@unine.ch

TÈl: +41 32 718 28 17; Fax: +41 32 718 28 01

\end{document}